\documentclass[twocolumn]{autart}    % Enable this line and disable the 
\usepackage{mathtools}
\usepackage[hidelinks]{hyperref}
\usepackage{algorithmic}
\usepackage{algorithm}
\usepackage{amsmath,amssymb}

\usepackage{color}
\usepackage{cite}
  
\allowdisplaybreaks[4]

\newtheorem{lemma}{Lemma}
\newtheorem{theorem}{Theorem}
\newtheorem{assumption}{Assumption}
\newtheorem{remark}{Remark}
\newtheorem{definition}{Definition}

\newtheorem{corollary}{Corollary}

\begin{document}

\begin{frontmatter}

\title{Projected Gradient Descent for Constrained
Decision-Dependent Optimization} % Title, preferably not more 

\thanks[footnoteinfo]{This work was supported in part by Swedish Research Council Distinguished Professor Grant 2017-01078, Knut and Alice Wallenberg Foundation, Wallenberg Scholar Grant, Swedish Strategic Research Foundation SUCCESS Grant FUS21-0026, AFOSR under award \#FA9550-19-1-0169, the Digital Futures Scholar-in-Residence Program, the European Union's Horizon 2020 research and innovation programme under grant agreement no. 739551 (KIOS CoE), and the Italian Ministry for Research in the framework of the 2017 Program for Research Projects of National Interest (PRIN), Grant no. 2017YKXYXJ.}

\author[DCS]{Zifan Wang}\ead{zifanw@kth.se},    % Add the 
\author[ECUST]{Changxin Liu}\ead{changxinl@ecust.edu.cn},               % e-mail address 
\author[IC]{Thomas Parisini}\ead{t.parisini@imperial.ac.uk},  % (ead) as shown
\author[DUKE]{Michael M. Zavlanos}\ead{michael.zavlanos@duke.edu},
\author[DCS]{Karl H. Johansson}\ead{kallej@kth.se}

\address[DCS]{Division of Decision and Control Systems, School of Electrical Engineering and Computer Science, KTH Royal Institute of Technology,  Stockholm, Sweden}    % Please supply                    
\address[ECUST]{The Key Laboratory of Smart Manufacturing in Energy Chemical Process, Ministry of Education, East China University of Science and Technology, Shanghai 200237, China}    % Please supply
\address[IC]{ Department of Electrical and Electronic Engineering, Imperial College London, London, UK}             % full addresses
\address[DUKE]{Mechanical Engineering and Material Science, Duke University, Durham, USA}       % here.

                             % keyword list or with the 
                                          % help of the Automatica 
                                          % keyword wizard

\begin{abstract}                          % Abstract of not more than 200 words.
This paper considers the decision-dependent optimization problem, where the data distributions react in response to decisions affecting both the objective function and linear constraints. We propose a new method termed repeated projected gradient descent (RPGD), which iteratively projects points onto evolving feasible sets throughout the optimization process. To analyze the impact of varying projection sets, we show a Lipschitz continuity property of projections onto varying sets with an explicitly given Lipschitz constant. Leveraging this property, we provide sufficient conditions for the convergence of RPGD to the constrained equilibrium point.  Compared to the existing dual ascent method, RPGD ensures continuous feasibility throughout the optimization process and reduces the computational burden. We validate our results through numerical experiments on a market problem and dynamic pricing problem.

\end{abstract}

\end{frontmatter}

\section{Introduction}

Machine learning models are typically optimized with the assumption that historically collected data will remain indicative of future system behavior. 
This assumption, however, often fails to hold in dynamic environments where the algorithm’s decisions can reshape the system it aims to model.
This phenomenon is observed in various applications ranging from 
online labor markets \cite{horton2010online} to vehicle-sharing systems  \cite{banerjee2015pricing,bianchin2021online,narang2022learning}, where strategic adjustments by users lead to shifts in data distribution.
To account for the distribution shift induced by the changes of decisions, recent advances have given rise to the frameworks of performative prediction \cite{perdomo2020performative,yan2024zero,miller2021outside} and optimization with decision-dependent distributions \cite{drusvyatskiy2023stochastic,wood2021online,wood2023stochastic}. 
These frameworks highlight the reciprocal nature of predictions and decision-making, where the predictive models can trigger changes in the outcomes they aim to predict.

Expectation constraints are prevalent in many machine learning problems, such as resource allocation \cite{lan2020algorithms,wang2023constrained}, fairness-constrained classification \cite{akhtar2021conservative,im2023stochastic} and financial risk management \cite{krokhmal2002portfolio}.
Works in these fields often operate under the assumption that the distribution of constraints remains static.
However, this assumption is not always true, as the distribution within the constraints can be influenced by the decisions made within these models.
An illustrative example is dynamic pricing for parking management.
In this scenario, the pricing decisions directly impact each user's parking time, which, in turn, can affect the system's constraints if the parking duration is used to define them.

Motivated by the discussions above, this paper considers decision-dependent distributions that occur in both the objective and linear constraints.
To the best of our knowledge, only \cite{wang2023constrained} solves this problem by proposing repeated constrained minimization (RCM) and repeated dual ascent (RDA) methods, both of which are supported by theoretical convergence guarantees.
However, each method has its drawbacks: RCM, while solving a constrained minimization problem exactly at each iteration, incurs significant computational overhead, and RDA, operating in the dual space, struggles to maintain the feasibility of solutions throughout the optimization process.
To address these limitations, we propose a repeated projected gradient descent (RPGD) method that iteratively updates decisions by projecting them onto evolving feasible sets. 
This method ensures the feasibility of solutions at each iteration thanks to the projection operator, although it results in time-varying feasible sets which may pose analytical challenges.
To understand the impact of these evolving projection sets, we begin by establishing the Lipschitz continuity of the projections. 
This foundational result allows us to derive sufficient conditions for the convergence of RPGD.
Notably, when the constraints are fixed, the conditions required for RPGD's convergence are less restrictive compared to those for RDA.
From a computational perspective, RPGD primarily involves computing gradients and projections each iteration, whereas RDA requires exact solutions of two unconstrained minimization problems. 
When the projection operations are straightforward, as is often the case with linear constraints, RPGD tends to be less computationally intensive.

% A closely related line of work explores constrained optimization with varying constraints using projected gradient descent methods \cite{wood2021online,cao2018online,yi2020distributed1}. 
% For example, \cite{wood2021online} addresses decision-dependent optimization with time-varying constraints and proposes a projected gradient descent method to track the evolving equilibrium point. 
% However, all the above works assume that the variability of constraints is dependent on time rather than the decisions. 
% Consequently, these methodologies fail to capture how decisions influence the evolution of constraints.
% This limitation renders their techniques not applicable in our case.
% We address this gap by exploring a Lipschitz continuity of projections defined by linear constraints.

The paper is organized as follows. In Section~\ref{sec:problem}, we formulate our problem and present the key assumptions. In Section~\ref{sec:alg}, we propose our algorithm and analyze its convergence. Section~\ref{sec:exp} provides the numerical validation of the proposed algorithm for a market problem and a dynamic pricing problem for parking management. Finally, we conclude this paper in Section~\ref{sec:conclusion}.

% Advantage:
% Guarantee feasibility; less restrict condition when $\epsilon_g$;

% projection QP; computation efficient when f is complex
\section{Problem Setup}\label{sec:problem}
This section presents the problem definition and key assumptions. 
Throughout the paper, we specify $\mathbb{R}^n$ as the $n$-dimensional Euclidean space equipped with an inner product $\langle \cdot, \cdot \rangle$ and a 2-norm $\left\|\cdot\right\|$. Furthermore, $\mathbb{R}^n$ is endowed with a $\sigma$-algebra. For a matrix $A \in \mathbb{R}^{m\times n}$, the notation $\left\| A\right\|_2$ represents the maximum singular value of $A$. For a positive definite matrix $Q\in \mathbb{R}^{n \times n}$, $\lambda_{\min}(Q)$ denotes the smallest eigenvalue of $Q$.

\subsection{Problem Definition}

We consider the linearly constrained optimization problem as in \cite{wang2023constrained}
\begin{align}\label{eq:problem}
    \min_x &\mathop{\mathbb{E}}_{z\sim \mathcal{D}(x)}[l(x,z)] \nonumber \\
    {\rm{s.t.}} & \quad Gx \leq \mathop{\mathbb{E}}_{w\sim \mathcal{D}_g(x)}[w],
\end{align}
where $x\in\mathbb{R}^n$ represents the decision variable, $G\in \mathbb{R}^{d_w\times n}$ is a matrix associated with the linear constraint.
The random variables $z\in \mathbb{R}^{d_z}$ and $w\in \mathbb{R}^{d_w}$ are associated with the objective function and constraints, respectively. 
The random variable $z$ is distributed according to $\mathcal{D}(x)$, a distribution map that transforms the decision variable space $\mathbb{R}^n$ into a distribution space, making the objective function dependent on $x$.
Furthermore, the constraints are also decision-dependent, with the random variable $w$ that depends on $x$.
We assume that the loss function $l: \mathbb{R}^{n} \times \mathbb{R}^{d_z} \rightarrow \mathbb{R}$ is twice continuously differentiable with respect to $(x,z)$.
Moreover, we assume that appropriate Borel measurability conditions are satisfied, ensuring that the expected value operators $\mathop{\mathbb{E}}\limits_{z\sim \mathcal{D}(x)}[\cdot]$ and $\mathop{\mathbb{E}}\limits_{z\sim \mathcal{D}_g(x)}[\cdot]$  are rigorously defined.

In this paper, we are interested in the following solution point.

\begin{definition}(Constrained Equilibrium Point).
\label{def:CSP}
A vector $x_s$ is called a constrained equilibrium point if it satisfies
\begin{align*}
    x_s =  \mathop{\rm{arg \; min}}_{x}  & \mathop{\mathbb{E}}_{z\sim \mathcal{D}(x_s)}[l(x,z)] \nonumber \\
     {\rm{s.t.}}& \quad Gx\leq \mathop{\mathbb{E}}_{w\sim \mathcal{D}_g(x_s)}[w] .
\end{align*}
\end{definition}

An equilibrium point $x_s$ is a fixed point that solves the constrained optimization problem using the distributions it induces. In the context of strategic classification, the decision variable $x_s$ is considered to be at equilibrium if the institution has no incentive to adopt a different classifier, based solely on the population's response to $x_s$. The equilibrium point is also referred to as a stable point in \cite{perdomo2020performative}.

Our goal in this paper is to design an algorithm that converges to the constrained equilibrium point while reducing computational overhead and maintaining feasibility during the optimization process.

\subsection{Key Assumptions}

The distribution maps in \eqref{eq:problem} are crucial in modeling how the distributions respond to changes in decisions. 
Convergence guarantees become unattainable when the distribution maps lack defined constraints.
Therefore, it is reasonable to make some regularity assumptions on these maps, commonly referred to as $\epsilon$-sensitivity.
\begin{definition}($\epsilon$-sensitivity).
We say that a distribution map $\mathcal{D}(\cdot)$ is $\epsilon$-sensitive if for all $x,y\in \mathbb{R}^n$, we have 
\begin{align}
    W_1(\mathcal{D}(x),\mathcal{D}(y))\leq \epsilon \left\|x-y\right\|,
\end{align}
where $W_1$ denotes the earth mover's distance.
\end{definition}
The earth mover's distance, also known as $1$-Wasserstein distance, is a measure of the minimum cost required to transform one distribution into another. 
We also make the following assumptions that are common in the decision-dependent optimization literature \cite{perdomo2020performative,miller2021outside,mendler2020stochastic}
\begin{assumption}\label{assump:strong_convex}
The loss function $l(x,z)$ is $\gamma$-strongly convex in $x$ for every $z \in \mathcal{Z}$, i.e., 
\begin{align*}
    l(x_1,z) \geq l(x_2,z) + \langle \nabla_x l(x_2,z),x_1-x_2 \rangle + \frac{\gamma}{2}\left\|x_1-x_2\right\|^2,
\end{align*}    
for all $x_1,x_2 \in \mathbb{R}^n$.
\end{assumption}

\begin{assumption}\label{assump:smooth}
$\nabla_x l(x,z)$ is $\beta_z$-Lipschitz continuous in $z$ for every $x\in \mathbb{R}^n$, i.e.,
\begin{align*}
    &\left\| \nabla_x l(x,z) - \nabla_x l(x,z')\right\| \leq \beta_z \left\| z-z'\right\|,
\end{align*}
for all $z,z'\in \mathcal{Z}$,
and $\nabla_x l(x,z)$ is $\beta_x$-Lipschitz continuous in $x$ for every $z\in \mathcal{Z}$, i.e.,
\begin{align*}
    &\left\| \nabla_x l(x_1,z) - \nabla_x l(x_2,z) \right\| \leq \beta_x \left\| x_1-x_2\right\|,
\end{align*}
for all $x_1,x_2\in \mathbb{R}^n$.
\end{assumption}

\begin{assumption} \label{assump:distribution_map}
The distribution maps $\mathcal{D}$ and $\mathcal{D}_g$ are $\epsilon$- and $\epsilon_g$-sensitive, respectively, namely,
\begin{align*}
    & W_1(\mathcal{D}(x),\mathcal{D}(y))\leq \epsilon \left\|x-y \right\|, \nonumber \\
    & W_1(\mathcal{D}_g(x),\mathcal{D}_g(y))\leq \epsilon_g \left\|x-y \right\|,
\end{align*}
for all $x,y \in \mathbb{R}^n$.
\end{assumption}
Moreover, we introduce the following assumption on the linear constraints as in \cite{wang2023constrained}.

\begin{assumption}\label{assumption:constraint}
    The matrix $G$ has full row rank.
\end{assumption}

Assumption \ref{assumption:constraint} is a technical condition that ensures linear convergence for constrained optimization problems, see e.g., \cite{wu2021new,qu2018exponential}.

\section{Algorithm}\label{sec:alg}

In this section, we propose a new algorithm that converges to the constrained equilibrium point. 

For simplicity of notation, we define functions
\begin{align*}
    f_{x'}(x) :=\mathop{\mathbb{E}}\limits_{z\sim \mathcal{D}(x')}[l(x,z)], \quad  \xi(x'):=\mathop{\mathbb{E}}\limits_{w\sim \mathcal{D}_g(x')}[w], 
\end{align*}
and the set 
\begin{align*}
    S(x')= \{ x\in \mathbb{R}^n: Gx \leq \xi(x')\}.
\end{align*}

\begin{algorithm}[t]
\caption{Repeated Projected Gradient Descent}
\begin{algorithmic}[1]
\STATE  \textbf{Input}: Initial variable $x_0$.
    \FOR {$t=0,1,2,\ldots$}
        \STATE $x_{t+1} = P_{S(x_t)} \left\{ x_t - \eta \nabla f_{x_t}(x_t)\right\}$;
    \ENDFOR
\STATE  \textbf{Output}: Sequence $\{ x_t\}$
\end{algorithmic}\label{alg:RPGD}
\end{algorithm}

The RPGD algorithm is presented in Algorithm~\ref{alg:RPGD}. At iteration $t+1$, projected gradient descent is applied to the constrained problem defined by the objective function $f_{x_t}$ and the feasible set $S(x_t)$, both of which are shaped by the distributions generated by the previous decision $x_t$.
Mathematically, RPGD performs the following update:
\begin{align}\label{eq:update}
    x_{t+1} = P_{S(x_t)} \left\{ x_t - \eta \nabla f_{x_t}(x_t)\right\},
\end{align}
where $\eta$ is the step size that will be determined later.
The projection operator of a vector $y\in\mathbb{R}^n$ onto the set $S$ is formally defined as 
\begin{align*}
    P_S\{y\} = \mathop{\rm{arg \; min}}\limits_{x\in \mathbb{R}^n} \delta_{C}(x) + \frac{1}{2} \left\|x - y \right\|^2,
\end{align*}
where the indicator function $\delta_{C}(x) = 0$ if $x\in C$ and $\infty$ otherwise.

One characteristic of the RPGD algorithm is that the projection set varies in response to changes in the decision variable. To understand the impact of these varying projection sets, we introduce the following lemma.

\begin{lemma}\label{lemma:proj_Lips}
Consider two closed, convex sets $C_1 = \{x: Gx \leq b_1 \}$ and $C_2 = \{x: Gx \leq b_2 \}$, where $G\in \mathbb{R}^{m\times n}$ has full row rank and $b_1, b_2 \in \mathbb{R}^{m}$ are constant vectors. Then, for any $y \in \mathbb{R}^{n}$, we have 
\begin{align}
    \left\| P_{C_1}(y) - P_{C_2}(y) \right\| \leq \frac{\left\| b_1 - b_2\right\|}{ \sqrt{\lambda_{\min}(G G^{\rm{T}})} }.
\end{align}
\end{lemma}

\textit{Proof of Lemma~\ref{lemma:proj_Lips}:}
Define $y_1 =P_{C_1}(y) $ and $y_2 =P_{C_2}(y)$.
By definition of projection, we have 
\begin{align}
    &y_1 = \mathop{\rm{arg \; min}}\limits_{x \in \mathbb{R}^n} \delta_{C_1}(x) + \frac{1}{2} \left\|x - y \right\|^2, \nonumber \\
    & y_2 = \mathop{\rm{arg \; min}}\limits_{x\in \mathbb{R}^n} \delta_{C_2}(x) + \frac{1}{2} \left\|x - y \right\|^2.
\end{align}
According to the first-order optimality condition, we have
\begin{align}\label{eq:proj:lemma:t1}
    & 0 \in \mathcal{N}_{C_1}(y_1) + y_1 -y, \nonumber \\
    & 0 \in \mathcal{N}_{C_2}(y_2) + y_2 -y,
\end{align}
where the normal cone $\mathcal{N}_{C}(y) = \{ z\in \mathbb{R}^n : \langle z,x-y \rangle \leq 0, \forall x\in C$. According to the definition of normal cone,  \eqref{eq:proj:lemma:t1} yields
\begin{align}
    &\langle y - y_1, x - y_1 \rangle \leq 0 , \forall x: Gx\leq b_1 \label{eq:proj:lemma:t2}, \\
    &\langle y - y_2, x - y_2 \rangle \leq 0 , \forall x: Gx\leq b_2 \label{eq:proj:lemma:t3}.
\end{align}
From the definition of projection, we have $y_1 \in C_1$, i.e., $G y_1 \leq b_1$. However, $y_2$ may not lie in the set $C_1$. For any point $h\in \mathbb{R}^n$ such that $G h = b_1 -b_2$, we have $G( y_2 + h) = G y_2 + b_1 - b_2 \leq b_1$. 
Therefore, the point $(y_2+h)$ lies in the set $C_1$ as long as $h$ satisfies $Gh = b_1 - b_2$.
Define the set $H=\{h\in \mathbb{R}^n: Gh = b_1 - b_2 \}$. Since the matrix $G$ has full rank, the set $H$ is always nonempty.  Similarly, $y_1$ may not lie in the set $C_2$. Using similar arguments, we can show that $(y_1-h)$ lies in the set $C_2$ for all $h\in H$. Substituting $(y_2 + h)$ and $(y_1-h)$ into $x$ in \eqref{eq:proj:lemma:t2} and \eqref{eq:proj:lemma:t3}, respectively, and adding these two inequalities together, we have 
\begin{align}\label{eq:proj:lemma:t4}
    \langle y_2 + h -y_1, y_2 - y_1 \rangle \leq 0,
\end{align}
for all $h\in H$. 
Rearranging \eqref{eq:proj:lemma:t4} yields 
\begin{align*}
    \left\| y_2 -y_1 \right\|^2 \leq \langle h, y_1 - y_2 \rangle \leq \left\| h\right\| \left\| y_1 -y_2 \right\|.
\end{align*}
Therefore, for all $h\in H$, we have
\begin{align}\label{eq:proj:lemma:t5}
    \left\| y_2 -y_1 \right\| \leq \left\|h\right\|.
\end{align}
We aim to find the tightest upper bound of $\left\|y_2 - y_1 \right\|$, which is equivalent to solving $\min\limits_{h\in H} \left\|h\right\|$. This problem can be reformulated as
\begin{align}\label{eq:proj:lemma:t6}
    \min_{h: Gh = b_1 - b_2} \frac{1}{2}\left\|h\right\|^2.
\end{align}
The problem \eqref{eq:proj:lemma:t6} is a quadratic programming problem with equality constraints and can be solved efficiently using the first-order necessary conditions. Let $h^*$ denote the optimal solution and $\lambda^*$ the associated Lagrange multiplier for \eqref{eq:proj:lemma:t6}. According to the KKT conditions, we have 
\begin{align}\label{eq:proj:lemma:t7}
    &h^* + G^{\rm{T}} \lambda^* =0, \nonumber \\
    &G x^* -(b_1-b_2)=0.
\end{align}
It is easy to verify that the solution of \eqref{eq:proj:lemma:t7} is unique and satisfies
\begin{align}
    &\lambda^* = - (GG^{\rm{T}})^{-1}(b_1 - b_2), \nonumber \\
    &h^* = -G^{\rm{T}} \lambda^* = G^{\rm{T}} (GG^{\rm{T}})^{-1}(b_1 - b_2).
\end{align}
Thus, the resulting optimal value satisfies 
\begin{align*}
    \left\| h^*\right\|^2 &= (b_1 - b_2)^{\rm{T}} ((GG^{\rm{T}})^{-1})^{\rm{T}} GG^{\rm{T}} (GG^{\rm{T}})-1 (b_1 - b_2) \\
    & = (b_1 - b_2)^{\rm{T}} (GG^{\rm{T}})^{-1} (b_1 - b_2) \\
    &\leq \left\| b_1 - b_2 \right\|^2 \left\| (GG^{\rm{T}})^{-1} \right\| \\
    &  = \left\| b_1 - b_2 \right\|^2 \frac{1}{\lambda_{\min}(G G^{\rm{T}})},
\end{align*}
where the last equality follows since $\left\|Q ^{-1}\right\| = \lambda_{\max}(Q^{-1}) =\frac{1}{\lambda_{\min}(Q)} $ for any positive definite matrix $Q$.

Substituting $h^*$ into \eqref{eq:proj:lemma:t5}, we have
\begin{align*}
    \left\| y_2 -y_1 \right\| \leq \frac{\left\| b_1 - b_2\right\|}{\sqrt{\lambda_{\min}(G G^{\rm{T}})}}.
\end{align*}
The proof is complete. \hfill $\qed$

Lemma~\ref{lemma:proj_Lips} shows that the projection of a vector onto two different sets satisfies a Lipschitz continuity property. 
According to \cite{bednarczuk2021lipschitz}, within the Hilbert space, the projection of any point $v$ onto two closed convex sets $C$ and $C'$ satisfies $\left\| P_{C}(v) - P_{C'}(v)\right\| \leq L_H d_{H}(C,C')$, where $d_{H}$ denotes the Hausdorff distance between two sets, see \cite{rockafellar2009variational} for more details about Hausdorff distance. 
The Lipschitz parameter $L_H$ depends on the specific space and the projection sets involved and is typically hard to determine.
In Lemma~\ref{lemma:proj_Lips}, we establish a specific instance of this Lipschitz property for linear constraints in Euclidian space and derive the explicit value of the Lipschitz constant.

Now we are ready to present the convergence result of RPGD. The proof can be found in the Appendix.

\begin{theorem}\label{theorem}
Suppose Assumptions~\ref{assump:strong_convex}--\ref{assumption:constraint} hold. If
\begin{align}\label{eq:RPGD:suff}
    \gamma - \epsilon \beta_z -  \frac{\epsilon_g (\epsilon \beta_z + \beta_x)}{\sqrt{\lambda_{\min}(G G^{\rm{T}})}} >0,
\end{align}
and 
\begin{align}\label{eq:RPGD:suff2}
    (\gamma - \epsilon \beta_z)^2 - \frac{2 \epsilon_g (\epsilon \beta_z + \beta_x) (\gamma+\beta_x)}{\sqrt{\lambda_{\min}(G G^{\rm{T}})}} >0,
\end{align}
then there exist $s_2>s_1\geq 0$ such that for all step size $\eta \in (s_1,s_2)$ we have 
\begin{align*}
    \left\| x_{t+1} - x_s \right\|^2 
    \leq \kappa^t \left\| x_1 - x_s \right\|^2,
\end{align*}
where 
\begin{align*}
    & \kappa = c_2 \eta^2 -2\eta c_1 + c_0 +1 < 1, \\
    &c_2 = ( \epsilon \beta_z + \beta_x)^2, \; c_1 = \gamma - \epsilon \beta_z - \frac{\epsilon_g (\epsilon \beta_z +\beta_x)}{\left\| G\right\|},  \nonumber \\
    &c_0 = \frac{\epsilon_g^2}{\left\| G\right\|^2} + \frac{2\epsilon_g}{\left\| G\right\|}, \; \Delta = 4c_1^2 - 4 c_2 c_0, \nonumber \\
    & s_1 = \frac{ 2c_1 - \sqrt{\Delta}}{2 c_2}, \; s_2 = \frac{ 2c_1 + \sqrt{\Delta}}{2 c_2}.
\end{align*}
\end{theorem}

Theorem~\ref{theorem} shows that RPGD converges to the constrained equilibrium point under conditions~\eqref{eq:RPGD:suff} and \eqref{eq:RPGD:suff2}. 
According to \cite{wang2023constrained}, when $\gamma - \epsilon \beta - \sqrt{\frac{\gamma \beta_x}{\lambda_{\rm{min} (GG^{\rm{T}})}}} \epsilon_g>0$, the constrained equilibrium point is unique and the
repeated constrained minimization will converge with the cost of solving an exact constrained optimization problem per iteration.
Note that the condition $\gamma - \epsilon \beta - \sqrt{\frac{\gamma \beta_x}{\lambda_{\rm{min} (GG^{\rm{T}})}}} \epsilon_g>0$ is less stringent than the condition  \eqref{eq:RPGD:suff}, because 
\begin{align*}
    &\gamma - \epsilon \beta_z -  \frac{\epsilon_g }{\sqrt{\lambda_{\min}(G G^{\rm{T}})}} (\epsilon \beta_z + \beta_x) \\
    \leq & \gamma - \epsilon \beta_z - \frac{\epsilon_g \beta_x}{\sqrt{\lambda_{\min}(G G^{\rm{T}})}} \\
    \leq & \gamma - \epsilon \beta_z -\sqrt{\frac{\gamma \beta_x}{\lambda_{\rm{min} (GG^{\rm{T}})}}} \epsilon_g,
\end{align*}
where the last inequality follows since $\gamma \leq \beta_x$. Therefore, the condition \eqref{eq:RPGD:suff} ensures the uniqueness of the equilibrium point.

% {\color{red} not clear to me here}

When $\epsilon_g = 0$, i.e., the constraints are fixed, the condition for RPGD to converge becomes $\gamma > \epsilon \beta_x$, which is less stringent than that required for the repeated dual ascent (RDA) method \cite{wang2023constrained}.
We discuss some other differences between RPGD and RDA in the following remark.

\begin{remark}
RPGD directly minimizes the objective function in the primal space, ensuring that constraints are consistently met and the solutions remain feasible throughout the optimization process.
In contrast, RDA performs gradient ascent in the dual space, which abstracts from the direct manipulation of primal variables. Although RDA is guaranteed to converge, the sequence of solutions do not necessarily satisfy the primal constraints.
Moreover, from a computational perspective, RDA requires solving two minimization problems per iteration while RPGD requires only a single gradient computation followed by a projection step. 
The projection onto linear constraints, which can be viewed as solving a quadratic programming problem, can be solved efficiently.
This enhances the computational efficiency of RPGD compared to RDA.
\end{remark}

When the objective function is fixed ($\epsilon=0$), implying that the distribution in the objective function is decision-independent, the convergence of RPGD is presented in the following result.

\begin{corollary}\label{corollary}
Suppose Assumptions~\ref{assump:strong_convex}--\ref{assumption:constraint} hold and $\epsilon=0$. If
\begin{align}\label{eq:RPGD:suff:cor}
    \epsilon_g < \frac{\gamma^2 \sqrt{\lambda_{\min}(G G^{\rm{T}})}} {2 \beta_x (\beta_x + \gamma)},
\end{align}
then there exist $\bar{s}_2>\bar{s}_1 \geq 0$ such that for all step size
$\eta \in (\bar{s}_1,\bar{s}_2)$  we have
\begin{align*}
    \left\| x_{t+1} - x_s \right\|^2 
    \leq \bar{\kappa}^t \left\| x_1 - x_s \right\|^2,
\end{align*}
where 
\begin{align*}
    &\bar{\kappa} = \bar{c}_2 \eta^2 -2\eta \bar{c}_1 + \bar{c}_0 +1 < 1, \\
    &\bar{c}_2 =  \beta_x^2, \; \bar{c}_1 = \gamma  - \frac{\epsilon_g \beta_x}{\sqrt{\lambda_{\min}(G G^{\rm{T}})}}, \nonumber \\
    &\bar{c}_0 = \frac{\epsilon_g^2}{\lambda_{\min}(G G^{\rm{T}})} + \frac{2\epsilon_g}{\sqrt{\lambda_{\min}(G G^{\rm{T}})}}, \nonumber \\
    & \bar{\Delta} = 4 \bar{c}_1^2 - 4 \bar{c}_2 \bar{c}_0,\; \bar{s}_1 = \frac{ 2\bar{c}_1 - \sqrt{\bar{\Delta}}}{2 \bar{c}_2}, \; \bar{s}_2 = \frac{ 2\bar{c}_1 + \sqrt{\bar{\Delta}}}{2 \bar{c}_2}.
\end{align*}
\end{corollary}
\textit{Proof of Corollary~\ref{corollary}:}
Note that \eqref{eq:RPGD:suff:cor} yields $\epsilon_g < \frac{\gamma^2 \sqrt{\lambda_{\min}(G G^{\rm{T}})}} {2 \beta_x (\beta_x + \gamma)} < \frac{\gamma \sqrt{\lambda_{\min}(G G^{\rm{T}})}}{\beta_x}$ and thus we have $\bar{c}_1>0$.
Following the steps in the proof of Theorem~\ref{theorem}, we can obtain the desired result. The detailed proof is omitted.  \hfill $\qed$

\begin{remark}
Corollary~\ref{corollary} shows that the selection of $\eta$ must satisfy $\bar{c}_2 \eta^2 -2\eta \bar{c}_1 + \bar{c}_0 <0$. 
From this inequality, it is evident that $\eta$ cannot be arbitrarily small given that $\bar{c}_0>0$, indicating $\epsilon_g >0$ and the constraints are decision-dependent.
The underlying reasoning can be elucidated as follows.
In the dynamics of RPGD, the decision update involves taking a step in the direction of the negative gradient, followed by a projection onto the feasible set that is evolving with the decision itself. 
When the step size of the decision update is relatively small compared to the magnitude of changes in the projection set, these alterations in the projection set may project the decision update towards the direction of gradient ascent rather than descent.
Therefore, ensuring that $\eta$ is adequately large is imperative to maintain the alignment of the projected decision with the gradient descent direction.
\end{remark}

\section{Numerical Experiments}\label{sec:exp}
\subsection{Market Problem}
We consider the market problem introduced in \cite{wang2023constrained} with the loss function $l(x,z)=-x_1 z_1 - x_2 z_2 $, where $x_i$ is the price of the $i$-th good and $z_i$ represents the random demand of the $i$-th good, for $i=1,2$.
The demand of the first good is defined as $z_1 = \zeta_1(x_1) - a_1 x_1 $, where the random variable $\zeta_1(x_1)$ follows a uniform distribution $U[\zeta_{L_1},\zeta_{R_1}+ \epsilon x_1]$. The distribution of the random demand for the first good depends on the price.
For the second good, the demand is modeled as $z_2 = \zeta_2 - a_2 x_2$ with $\zeta_2$ following a constant uniform distribution $U[\zeta_{L_2}, \zeta_{R_2}]$.
The production costs $v_i$ for each good are also considered. The cost $v_1$ of the first good follows a uniform distribution $U[\underline{v}_1, 1.2 \underline{v}_1 + \epsilon_g x_1]$, implying that lower prices potentially lead to higher demands and thus, a lower average production cost.
The cost $v_2$ of the second good is subject to a uniform distribution $U[ \underline{v}_2, 1.2 \underline{v}_2]$.
The constraint is introduced as $Gx\leq \xi(x)$, where $G=[-a_3, -a_4]$, $\xi(x)=\mathbb{E}[-v_1(x_1)-v_2 -e_1]$. This constraint ensures that the selling prices exceed the weighted average production costs plus an additional expense $e_1$.

\begin{figure}[t]
\begin{center}
\centerline{\includegraphics[width=0.95\columnwidth]{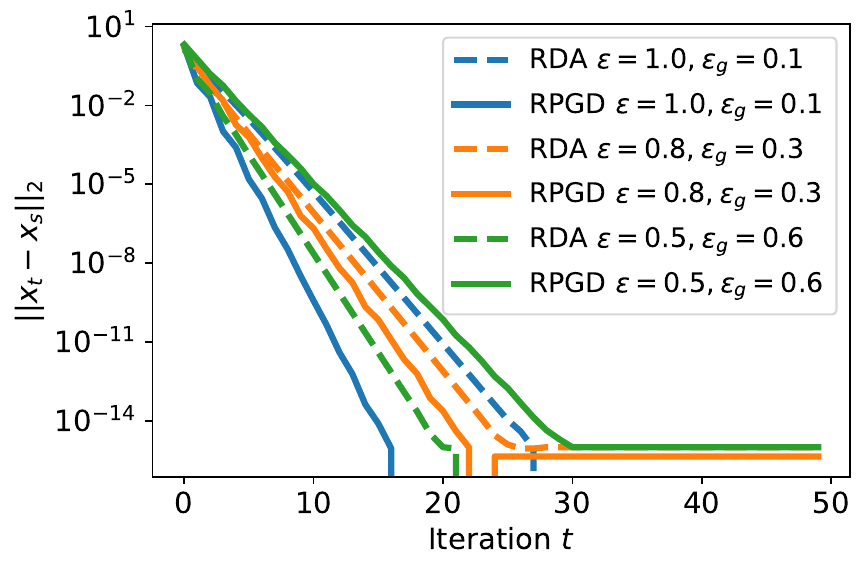}}
\caption{Convergence of RPGD and RDA for different values of $\epsilon$ and $\epsilon_g$ in the market problem.
}
\label{RPGD:market}
\end{center}
\end{figure}

\begin{figure}[t]
\begin{center}
\centerline{\includegraphics[width=0.95\columnwidth]{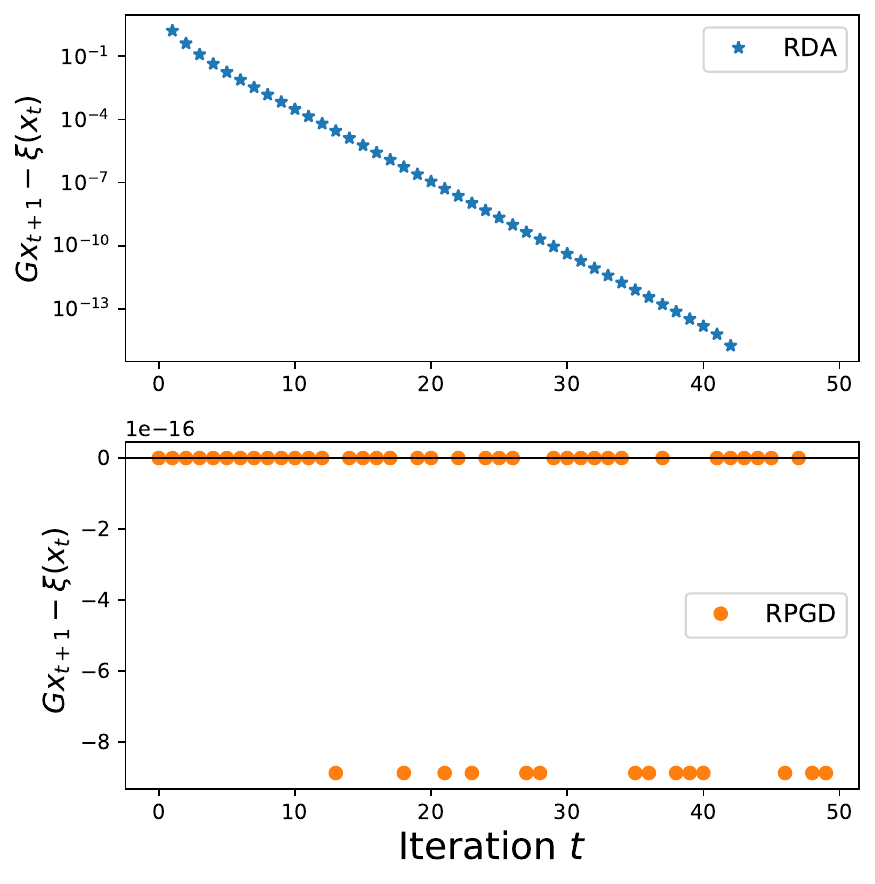}}
\caption{Constraint violation of RPGD and RDA when $\epsilon=1.5$, $\epsilon_g=0.4$ in the market problem. The horizontal black line represents $y=0$. The circle point denotes a feasible point while the star point denotes an infeasible point.} 
% {\color{red} why use dots for the constraint?}
\label{RPGD:constraint:vio}
\end{center}
\end{figure}

We set the following parameter values: $a_1=0.8$, $a_2=0.2$, $a_3=0.6$, $a_4 =1$, $\zeta_{L_1}=1$, $\zeta_{R_1}=5.5$, $\zeta_{L_2}=0.5$, $\zeta_{R_2}=2.2$, $e_1=1.2$, $\underline{v}_1=1.7$, $\underline{v}_2=2.5$. 
We compare the convergence of RPGD with RDA \cite{wang2023constrained} for different values of $\epsilon$ and $\epsilon_g$ and the results are displayed in Fig.~\ref{RPGD:market}. We observe that RPGD demonstrates faster convergence than RDA when $\epsilon_g$ is small, indicating that smaller changes in the constraints significantly enhance the efficiency of RPGD.
We investigate the constraint violation of decision updates during the optimization process by examining the values of $Gx_{t+1} - \xi(x_t)$, where negative values indicate feasible solutions.
We set $\epsilon=1.5$ and $\epsilon_g=0.4$. 
As shown in Fig~\ref{RPGD:constraint:vio}, RPGD consistently maintains the feasibility of decision updates at each iteration, whereas RDA fails at most iterations.

\subsection{Dynamic Pricing Problem}
We consider the application of a dynamic pricing experiment based on the model presented in \cite{wang2023constrained} and using the dataset from \cite{sfpark_data}.
The loss function is formulated as
$l(x,z_1,z_2)= (z_1 - 0.7)^2 - t z_2(x+ \bar{x}) + \frac{v}{2} \left\|x\right\|^2$, where $x$ is the price adjustment from the base price $\bar{x}=3$, $z_1$ denotes the occupancy rate, $z_2$ represents the total parking time. 
The goal is to maintain occupancy around 70\% and simultaneously maximize the parking revenue $z_2 (x+\bar{x})$. 
The occupancy rate is modeled by $z_1 = \zeta-Ax$, where $A=0.157$ and $\zeta $ drawn from a fixed distribution $P_0$ derivaed from the dataset.
When analyzing the individual response to a price adjustment $x$, we assume that each user adapts their behavior according to a best response strategy, resulting in the updated parking time $z_2' = z_2 - \epsilon x$. 
Additionally, we impose a constraint on the expected parking time as follows: $\mathop{\mathbb{E}}\limits_{z_2 \sim \mathcal{D}(x)}[z_2] \leq c_1 x +c_2$, meaning that the total occupied time should be linearly bounded. This constraint aims to balance economic efficiency with practical feasibility in parking management.

\begin{figure}[t]
\begin{center}
\centerline{\includegraphics[width=0.95\columnwidth]{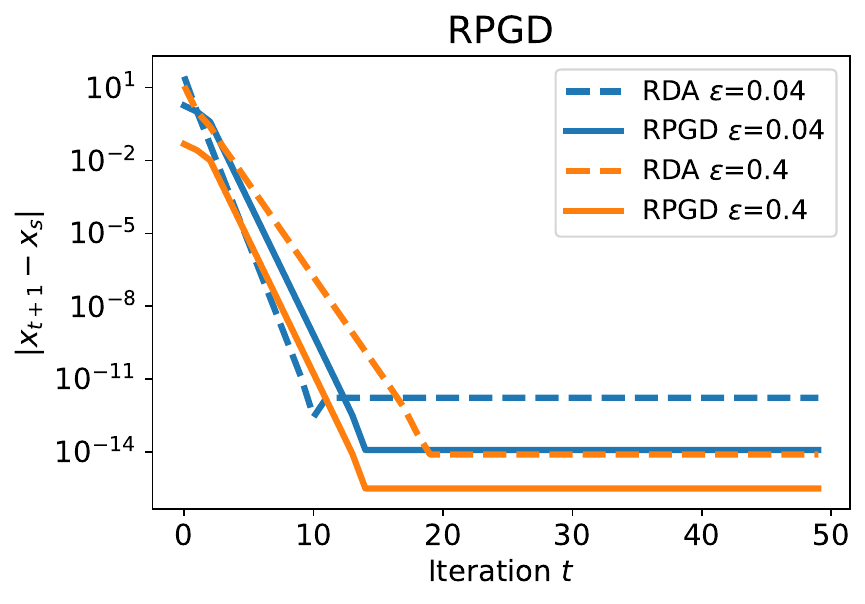}}
\caption{Convergence of RPGD and RDA for different values of $\epsilon$ and $\epsilon_g$ in the dynamic parking pricing problem.
}
\label{RPGD:park}
\end{center}
\end{figure}

\begin{figure}[t]
\begin{center}
\centerline{\includegraphics[width=0.95\columnwidth]{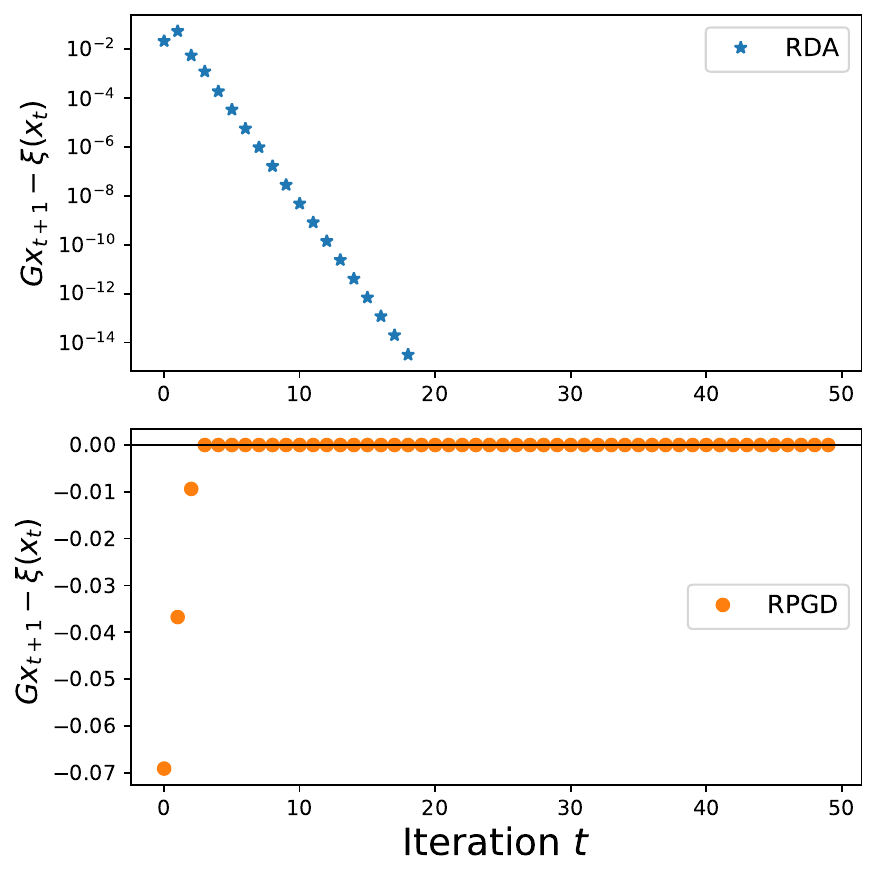}}
\caption{Constraint violation of RPGD and RDA when $\epsilon=0.4$ in the dynamic parking pricing problem. The horizontal black line represents $y=0$. The circle point denotes a feasible point while the star point denotes an infeasible point. }
\label{RPGD:park:constraint:vio}
\end{center}
\end{figure}

The parameters for the optimization process are set as follows: $v=0.03$, $t=0.005$, $c_1 =0.5$, $c_2 = 5$, with an initial price adjustment $x_0=0$. The step size for each algorithm is optimally tuned. 
Figure~\ref{RPGD:park} presents the convergence results for RPGD and RDA, demonstrating that both algorithms achieve convergence within a limited number of iterations. However, as shown in Fig.~\ref{RPGD:park:constraint:vio}, while RPGD maintains feasibility throughout the optimization process, RDA fails to meet the feasibility constraints during the initial 20 iterations.

\section{Conclusion}\label{sec:conclusion}
This paper proposed the RPGD method to address  optimization problems with decision-dependent constraints.
Although the projections sets are varying during the optimization process, we still showed the convergence of RPGD based on the constructed Lipschitz continuity property of projection onto varying sets.
The effectiveness of RPGD was validated through numerical experiments on market and dynamic pricing problems, which confirmed its superior performance in maintaining feasibility. 
Future work may explore the extension of RPGD to non-linear constraints and its application to broader classes of decision-dependent problems.
\section*{Appendix. Proof of Theorem~\ref{theorem}}\label{sec:appendix}
We first present a lemma that helps the analysis.
\begin{lemma}\label{lemma:xs}
For a constrained equilibrium point $x_s$, it satisfies 
\begin{align}
    x_s = P_{S(x_s)}\{ x_s -\eta \nabla f_{x_s}(x_s)\}  ,
\end{align}
for all $\eta \geq 0$.
\end{lemma}
\textit{Proof of Lemma~\ref{lemma:xs}:}
Assume $y = P_{S(x_s)}\{ x_s -\eta \nabla f_{x_s}(x_s)\}$. According to the definition of projection, we have 
\begin{align*}
    y = \mathop{\rm{arg \; min}}\limits_{x} \delta_{S(x_s)}(x) + \frac{1}{2} \left\|x - x_s + \eta \nabla f_{x_s}(x_s) \right\|^2,
\end{align*}
which is equivalent to 
\begin{align*}
    0 \in \mathcal{N}_{S(x_s)}(y) + y - x_s + \eta\nabla f_{x_s}(x_s).
\end{align*}
By virtue of the definition of the normal cone, we have
\begin{align}\label{eq:lemma:xs:t1}
    \langle x_s - \eta \nabla f_{x_s}(x_s) - y,x-y\rangle \leq 0, \forall x \in S(x_s).
\end{align}
Since $x_s \in S(x_s)$, substituting $x$ with $x_s$ in \eqref{eq:lemma:xs:t1} yields 
\begin{align}\label{eq:lemma:xs:t2}
    0\leq \left\| x_s - y \right\|^2 \leq \eta \langle \nabla f_{x_s}(x_s),x_s -y \rangle.
\end{align}
Due to the first-order optimality condition, we have 
\begin{align}\label{eq:lemma:xs:t3}
    \langle \nabla f_{x_s}(x_s), x - x_s \rangle \geq 0 , \forall x \in S(x_s).
\end{align}
Since $y\in S(x_s)$, substituting $x$ with $y$ in \eqref{eq:lemma:xs:t3} yields
\begin{align}\label{eq:lemma:xs:t4}
    \langle \nabla f_{x_s}(x_s), y - x_s \rangle \geq 0.
\end{align}
Combining \eqref{eq:lemma:xs:t2} and \eqref{eq:lemma:xs:t4}, we have
\begin{align*}
    0\leq \left\| x_s - y \right\|^2 \leq \eta \langle \nabla f_{x_s}(x_s),x_s -y \rangle \leq 0.
\end{align*}
Therefore, every inequality becomes equality, which yields $y=x_s$. The proof is complete.   \hfill $\qed$

By the update rule \eqref{eq:update}, we have
\begin{align}\label{eq:Thm:pf:t0}
    &\left\| x_{t+1} - x_s \right\|^2  = \left\| P_{S(x_t)} \left\{ x_t - \eta \nabla f_{x_t}(x_t)\right\} - x_s \right\|^2 
    \nonumber \\
    & = \left\| P_{S(x_t)} \left\{ x_t - \eta \nabla f_{x_t}(x_t)\right\} - P_{S(x_s)} \left\{ x_t - \eta \nabla f_{x_t}(x_t)\right\} \right. \nonumber \\
    &\quad \left. + P_{S(x_s)} \left\{ x_t - \eta \nabla f_{x_t}(x_t)\right\}- x_s \right\|^2  \nonumber \\
    & = \left\| P_{S(x_t)} \left\{ x_t - \eta \nabla f_{x_t}(x_t)\right\} - P_{S(x_s)} \left\{ x_t - \eta \nabla f_{x_t}(x_t)\right\} \right\|^2  \nonumber \\
    &\quad + 2 \Big\langle P_{S(x_s)} \left\{ x_t - \eta \nabla f_{x_t}(x_t)\right\}- x_s, \nonumber \\
    &\quad P_{S(x_t)} \left\{ x_t - \eta \nabla f_{x_t}(x_t)\right\} - P_{S(x_s)} \left\{ x_t - \eta \nabla f_{x_t}(x_t)\right\}  \Big\rangle \nonumber \\
    &\quad + \left\| P_{S(x_s)} \left\{ x_t - \eta \nabla f_{x_t}(x_t)\right\}- x_s \right\|^2 .
\end{align}
By virtue of Lemma~\ref{lemma:proj_Lips}, we have 
\begin{align}\label{eq:Thm:pf:t1}
    &\left\| P_{S(x_t)} \left\{ x_t - \eta \nabla f_{x_t}(x_t)\right\} - P_{S(x_s)} \left\{ x_t - \eta \nabla f_{x_t}(x_t)\right\} \right\| \nonumber \\
    &\leq \frac{ \left\| \xi(x_t) - \xi(x_s) \right\|} { \sqrt{\lambda_{\min}(G G^{\rm{T}})}} \leq \epsilon_g \frac{ \left\| x_t - x_s \right\|} { \sqrt{\lambda_{\min}(G G^{\rm{T}})} }.
\end{align}
Note that 
\begin{align}\label{eq:Thm:pf:t2}
    &\left\|  \nabla f_{x_t}(x_t) -  \nabla f_{x_s}(x_s) \right\| \nonumber \\
    & = \left\|  \nabla f_{x_t}(x_t) - \nabla f_{x_s}(x_t) +  \nabla f_{x_s}(x_t)- \nabla f_{x_s}(x_s) \right\| \nonumber \\
    &\leq \left\|  \nabla f_{x_t}(x_t) - \nabla f_{x_s}(x_t) \right\| + \left\| \nabla f_{x_s}(x_t)- \nabla f_{x_s}(x_s) \right\| \nonumber \\
    &\leq (\epsilon \beta_z + \beta_x) \left\| x_t - x_s\right\|,
\end{align}
where the last inequality follows from the smooth assumption and Lemma 2.1 in \cite{drusvyatskiy2023stochastic}.
By virtue of Lemma~\ref{lemma:xs}, we have 
\begin{align}\label{eq:Thm:pf:t3}
    &2\Big\langle P_{S(x_t)} \left\{ x_t - \eta \nabla f_{x_t}(x_t)\right\} - P_{S(x_s)} \left\{ x_t - \eta \nabla f_{x_t}(x_t)\right\}, \nonumber \\
    &\quad P_{S(x_s)} \left\{ x_t - \eta \nabla f_{x_t}(x_t)\right\}- x_s  \Big\rangle \nonumber \\ 
    & = 2\Big\langle P_{S(x_t)} \left\{ x_t - \eta \nabla f_{x_t}(x_t)\right\} - P_{S(x_s)} \left\{ x_t - \eta \nabla f_{x_t}(x_t)\right\}, \nonumber \\
    &\quad P_{S(x_s)} \left\{ x_t - \eta \nabla f_{x_t}(x_t)\right\}- P_{S(x_s)}\{ x_s -\eta \nabla f_{x_s}(x_s)\}  \Big\rangle \nonumber \\ 
    &\leq 2 \left\| P_{S(x_t)} \left\{ x_t - \eta \nabla f_{x_t}(x_t)\right\} - P_{S(x_s)} \left\{ x_t - \eta \nabla f_{x_t}(x_t)\right\} \right\| \nonumber \\
    & \quad \left\| P_{S(x_s)} \left\{ x_t - \eta \nabla f_{x_t}(x_t)\right\}- P_{S(x_s)}\{ x_s -\eta \nabla f_{x_s}(x_s)\}\right\| \nonumber \\
    & \leq  2 \left\| P_{S(x_t)} \left\{ x_t - \eta \nabla f_{x_t}(x_t)\right\} - P_{S(x_s)} \left\{ x_t - \eta \nabla f_{x_t}(x_t)\right\} \right\| \nonumber \\
    &\quad  \left\|  x_t - \eta \nabla f_{x_t}(x_t) -  x_s +\eta \nabla f_{x_s}(x_s) \right\| \nonumber \\ 
    &\leq  \frac{ 2\epsilon_g\left\| x_t - x_s \right\|} { \sqrt{\lambda_{\min}(G G^{\rm{T}})} } \left\|  x_t - \eta \nabla f_{x_t}(x_t) -  x_s +\eta \nabla f_{x_s}(x_s) \right\| \nonumber \\
    & \leq  \frac{ 2\epsilon_g \left\| x_t - x_s \right\|} { \sqrt{\lambda_{\min}(G G^{\rm{T}})} } \Big( \left\|  x_t - x_s \right\| + \eta \left\|  \nabla f_{x_t}(x_t) - \nabla f_{x_s}(x_s) \right\| \Big) \nonumber \\
    &\leq   \frac{2\epsilon_g  } { \sqrt{\lambda_{\min}(G G^{\rm{T}})} } \Big(1+ \eta(\epsilon \beta_z + \beta_x) \Big)\left\| x_t - x_s \right\|^2,
\end{align}
where the second inequality holds since the projection operator is non-expansive and the third inequality follows from \eqref{eq:Thm:pf:t1}. The last inequality follows from \eqref{eq:Thm:pf:t2}.
Moreover, we have
\begin{align}\label{eq:Thm:pf:t4}
    &\left\| P_{S(x_s)} \left\{ x_t - \eta \nabla f_{x_t}(x_t)\right\}- x_s \right\|^2  \nonumber \\
    & = \left\| P_{S(x_s)} \left\{ x_t - \eta \nabla f_{x_t}(x_t)\right\}- P_{S(x_s)}\{ x_s -\eta \nabla f_{x_s}(x_s) \}\right\|^2 \nonumber \\
    & \leq \left\| x_t - \eta \nabla f_{x_t}(x_t) - x_s +\eta \nabla f_{x_s}(x_s) \right\|^2 \nonumber \\
    & = \left\| x_t - x_s \right\|^2  + \eta^2 \left\| \nabla f_{x_t}(x_t) - \nabla f_{x_s}(x_s) \right\|^2 \nonumber \\
    &\quad - 2\eta \langle x_t-x_s, \nabla f_{x_t}(x_t) - \nabla f_{x_s}(x_s) \rangle  \nonumber \\
    & \leq \left\| x_t - x_s \right\|^2 + (\epsilon \beta_z + \beta_x)^2 \eta^2 \left\| x_t - x_s\right\|^2 \nonumber \\
    &\quad - 2\eta \langle x_t-x_s, \nabla f_{x_t}(x_t) - \nabla f_{x_s}(x_t) \rangle  \nonumber \\
    &\quad - 2 \eta \langle x_t-x_s,\nabla f_{x_s}(x_t) - \nabla f_{x_s}(x_s) \rangle  \nonumber \\
    &\leq \left\| x_t - x_s \right\|^2 + \eta^2(\epsilon \beta_z + \beta_x)^2  \left\| x_t - x_s\right\|^2  \nonumber \\
    &\quad +  2 \eta \epsilon \beta_z\left\| x_t - x_s\right\|^2 - 2\eta \gamma \left\| x_t - x_s\right\|^2 \nonumber \\
    &=\big( 1+ \eta^2(\epsilon \beta_z + \beta_x)^2 + 2 \eta \epsilon \beta_z - 2\eta \gamma \big) \left\| x_t - x_s\right\|^2.
\end{align}
The first inequality holds since the projection operator is non-expansive, the second inequality follows from \eqref{eq:Thm:pf:t2} and the third inequality follows from strong convexity and \eqref{eq:Thm:pf:t2}.

Substituting \eqref{eq:Thm:pf:t1}, \eqref{eq:Thm:pf:t3} and \eqref{eq:Thm:pf:t4} into \eqref{eq:Thm:pf:t0}, we have
\begin{align}\label{eq:Thm:pf:t5}
    &\left\| x_{t+1} - x_s \right\|^2 \nonumber \\
    &\leq  \frac{ \epsilon_g^2} { \lambda_{\min}(G G^{\rm{T}})}   \left\| x_t - x_s \right\|^2 \nonumber \\
    &\quad + \frac{2\epsilon_g  } { \sqrt{\lambda_{\min}(G G^{\rm{T}})} } \Big(1+ \eta(\epsilon \beta_z + \beta_x) \Big) \left\| x_t - x_s \right\|^2 \nonumber \\
    &\quad +\big( 1+ \eta^2(\epsilon \beta_z + \beta_x)^2 + 2 \eta \epsilon \beta_z - 2\eta \gamma \big) \left\| x_t - x_s\right\|^2 \nonumber \\
    & = \Big( \eta ^2 \big( \epsilon \beta_z + \beta_x \big)^2  -2\eta \big( \gamma - \epsilon \beta_z - \frac{\epsilon_g (\epsilon \beta_z +\beta_x)}{\sqrt{\lambda_{\min}(G G^{\rm{T}})}} \big) \nonumber \\
    &\quad + \big(\frac{\epsilon_g}{\sqrt{\lambda_{\min}(G G^{\rm{T}})} } + 1 \big)^2 \Big)\left\| x_t - x_s \right\|^2 \nonumber \\
    & = ( c_2 \eta^2 - 2 \eta c_1 + c_0 +1)\left\| x_t - x_s \right\|^2,
\end{align}
where we define $c_2 = ( \epsilon \beta_z + \beta_x)^2$, $c_1 = \gamma - \epsilon \beta_z - \frac{\epsilon_g (\epsilon \beta_z +\beta_x)}{\sqrt{\lambda_{\min}(G G^{\rm{T}})}}$, $c_0 = \frac{\epsilon_g^2}{\lambda_{\min}(G G^{\rm{T}})} + \frac{2\epsilon_g}{\sqrt{\lambda_{\min}(G G^{\rm{T}})}}$.
To show the convergence, it suffices to show that the choice of $\eta$ assures $c_2 \eta^2 -2\eta c_1 + c_0<0$.  Note that \eqref{eq:RPGD:suff} guarantees $c_1>0$. The discriminant of this quadratic function of $\eta$ satisfies
\begin{align*}
    &\Delta = 4c_1^2 - 4 c_2 c_0  \nonumber \\
    & = 4 \big( \gamma - \epsilon \beta_z - \frac{\epsilon_g (\epsilon \beta_z +\beta_x)}{ \sqrt{\lambda_{\min}(G G^{\rm{T}})}}\big)^2 \nonumber \\
    & \quad - 4( \epsilon \beta_z + \beta_x)^2 \big( \frac{\epsilon_g^2}{\lambda_{\min}(G G^{\rm{T}})} + \frac{2\epsilon_g}{\sqrt{\lambda_{\min}(G G^{\rm{T}})}}\big)  \nonumber \\
    & = 4\Big( (\gamma - \epsilon \beta_z)^2 - \frac{2\epsilon_g}{\sqrt{\lambda_{\min}(G G^{\rm{T}})}}(\epsilon \beta_z + \beta_x)(\gamma +\beta_x) \Big) \nonumber \\
    &>0,
\end{align*}
where the inequality follows from \eqref{eq:RPGD:suff2}. Since $\Delta>0$, there exist two distinct solution to $c_2 \eta^2 -2\eta c_1 + c_0 =0$, which we denote by $s_1,s_2$ with $s_1 = \frac{ 2c_1 - \sqrt{\Delta}}{2 c_2}$ and $s_2 = \frac{ 2c_1 + \sqrt{\Delta}}{2 c_2}$. When $\eta \in (s_1,s_2)$, we have $\kappa = c_2 \eta^2 -2\eta c_1 + c_0 +1 < 1 $. From \eqref{eq:Thm:pf:t5}, we have
\begin{align}
    \left\| x_{t+1} - x_s \right\|^2 \leq \kappa \left\| x_t -x_s \right\|^2 
    \leq \kappa^t \left\| x_1 - x_s \right\|^2.
\end{align}
The proof is complete.  \hfill $\qed$

\bibliography{autosam}

\begin{thebibliography}{10}

\bibitem{horton2010online}
John~J Horton.
\newblock Online labor markets.
\newblock In {\em International Workshop on Internet and Network Economics}, pages 515--522, 2010.

\bibitem{banerjee2015pricing}
Siddhartha Banerjee, Carlos Riquelme, and Ramesh Johari.
\newblock Pricing in ride-share platforms: A queueing-theoretic approach.
\newblock In {\em ACM Conference on Economics and Computation}, 2015.

\bibitem{bianchin2021online}
Gianluca Bianchin, Miguel Vaquero, Jorge Cortes, and Emiliano Dall'Anese.
\newblock Online stochastic optimization for unknown linear systems: Data-driven synthesis and controller analysis.
\newblock {\em arXiv preprint arXiv:2108.13040}, 2021.

\bibitem{narang2022learning}
Adhyyan Narang, Evan Faulkner, Dmitriy Drusvyatskiy, Maryam Fazel, and Lillian Ratliff.
\newblock Learning in stochastic monotone games with decision-dependent data.
\newblock In {\em International Conference on Artificial Intelligence and Statistics}, pages 5891--5912. PMLR, 2022.

\bibitem{perdomo2020performative}
Juan Perdomo, Tijana Zrnic, Celestine Mendler-D{\"u}nner, and Moritz Hardt.
\newblock Performative prediction.
\newblock In {\em International Conference on Machine Learning}, pages 7599--7609. PMLR, 2020.

\bibitem{yan2024zero}
Wenjing Yan and Xuanyu Cao.
\newblock Zero-regret performative prediction under inequality constraints.
\newblock In {\em Advances in Neural Information Processing Systems}. PMLR, 2024.

\bibitem{miller2021outside}
John~P Miller, Juan~C Perdomo, and Tijana Zrnic.
\newblock Outside the echo chamber: Optimizing the performative risk.
\newblock In {\em International Conference on Machine Learning}, pages 7710--7720. PMLR, 2021.

\bibitem{drusvyatskiy2023stochastic}
Dmitriy Drusvyatskiy and Lin Xiao.
\newblock Stochastic optimization with decision-dependent distributions.
\newblock {\em Mathematics of Operations Research}, 48(2):954--998, 2023.

\bibitem{wood2021online}
Killian Wood, Gianluca Bianchin, and Emiliano Dall’Anese.
\newblock Online projected gradient descent for stochastic optimization with decision-dependent distributions.
\newblock {\em IEEE Control Systems Letters}, 6:1646--1651, 2021.

\bibitem{wood2023stochastic}
Killian Wood and Emiliano Dall’Anese.
\newblock Stochastic saddle point problems with decision-dependent distributions.
\newblock {\em SIAM Journal on Optimization}, 33(3):1943--1967, 2023.

\bibitem{lan2020algorithms}
Guanghui Lan and Zhiqiang Zhou.
\newblock Algorithms for stochastic optimization with function or expectation constraints.
\newblock {\em Computational Optimization and Applications}, 76(2):461--498, 2020.

\bibitem{wang2023constrained}
Zifan Wang, Changxin Liu, Thomas Parisini, Michael~M Zavlanos, and Karl~H Johansson.
\newblock Constrained optimization with decision-dependent distributions.
\newblock {\em arXiv preprint arXiv:2310.02384}, 2023.

\bibitem{akhtar2021conservative}
Zeeshan Akhtar, Amrit~Singh Bedi, and Ketan Rajawat.
\newblock Conservative stochastic optimization with expectation constraints.
\newblock {\em IEEE Transactions on Signal Processing}, 69:3190--3205, 2021.

\bibitem{im2023stochastic}
Hyungki Im and Paul Grigas.
\newblock Stochastic first-order algorithms for constrained distributionally robust optimization.
\newblock {\em arXiv preprint arXiv:2305.16584}, 2023.

\bibitem{krokhmal2002portfolio}
Pavlo Krokhmal, Jonas Palmquist, and Stanislav Uryasev.
\newblock Portfolio optimization with conditional value-at-risk objective and constraints.
\newblock {\em Journal of Risk}, 4:43--68, 2002.

\bibitem{mendler2020stochastic}
Celestine Mendler-D{\"u}nner, Juan Perdomo, Tijana Zrnic, and Moritz Hardt.
\newblock Stochastic optimization for performative prediction.
\newblock In {\em Advances in Neural Information Processing Systems}, volume~33, pages 4929--4939, 2020.

\bibitem{wu2021new}
Xuyang Wu, Sindri Magn{\'u}sson, and Mikael Johansson.
\newblock A new family of feasible methods for distributed resource allocation.
\newblock In {\em 2021 60th IEEE Conference on Decision and Control}, pages 3355--3360. IEEE, 2021.

\bibitem{qu2018exponential}
Guannan Qu and Na~Li.
\newblock On the exponential stability of primal-dual gradient dynamics.
\newblock {\em IEEE Control Systems Letters}, 3(1):43--48, 2018.

\bibitem{bednarczuk2021lipschitz}
Ewa~M Bednarczuk and Krzysztof~E Rutkowski.
\newblock On lipschitz continuity of projections onto polyhedral moving sets.
\newblock {\em Applied Mathematics \& Optimization}, 84(2):2147--2175, 2021.

\bibitem{rockafellar2009variational}
R~Tyrrell Rockafellar and Roger J-B Wets.
\newblock {\em Variational Analysis}, volume 317.
\newblock Springer Science \& Business Media, 2009.

\bibitem{sfpark_data}
Sfpark parking sensor data hourly occupancy 2011 -- 2013.
\newblock {\em https://www.sfmta.com/getting-around/drive-park/demand-responsive-pricing/sfpark-evaluation}, 2013.

\end{thebibliography}
\bibliographystyle{unsrt} 

\end{document}